\documentclass [12pt] {article}
\usepackage{amsfonts}
\usepackage{amssymb}

\newcommand{\qed} {\hspace {0.1in} \rule {1.5mm} {3.5mm}}

\newtheorem{lemma}{Lemma}[section]
\newtheorem{corollary}{Corollary}[section]
\newtheorem{theorem}{Theorem}
\newtheorem{statement}{Statement}

\newtheorem{proposition}{Proposition}[section]

\def\limn{\lim_{n\to\infty}}

\def\e{\epsilon}

\def\ke{\mbox{Ker}\,}

\def\tr{\mbox{Tr}\,}
\def\mat{\mbox{Mat}\,}

\def\dim{{\rm dim}}
\def\supp{{\rm supp}}
\def\Fo{F$\mbox{\o}$lner}
\def\lan{\langle}
\def\ran{\rangle}

\def\i{^\infty_{n=1}}
\def\proof{\smallskip\noindent{\it Proof.} }

\def\bR{{\mathbb R}}
\def\bN{{\mathbb N}}
\def\bC{{\mathbb C}}
\def\bQ{{\mathbb Q}}

\def\cA{\mbox{$\cal A$}}

\def\cN{\mbox{$\cal N$}}

\def\cN{\mbox{$\cal N$}}

\def\lan{\langle}

\def\ran{\rangle}

\def\to{\rightarrow}

\def\osb{\overline{\sigma_B}}
\def\usb{\underline{\sigma_B}}
\begin{document}
\title{$L^2$-spectral invariants and quasi-crystal graphs}
\author{\sc G\'abor Elek
\footnote {The Alfred Renyi Mathematical Institute of
the Hungarian Academy of Sciences, P.O. Box 127, H-1364 Budapest, Hungary.
email:elek@renyi.hu, Supported by OTKA Grants T 049841 and T 037846}}
\date{}
\maketitle \vskip 0.2in \noindent{\bf Abstract.} Introducing and
studying the pattern frequency algebra, we  prove the analogue of
L\"uck's approximation theorems on $L^2$-spectral invariants in the
case of aperiodic order. These results imply a uniform convergence
theorem for the integrated density of states as well as the
positivity of the logarithmic determinant of
certain discrete Schrodinger operators.

\vskip 0.2in
\noindent{\bf AMS Subject Classifications:} 81Q10, 46L51
\vskip 0.2in
\noindent{\bf Keywords:\,}quasi-crystal graphs, $L^2$-invariants, von Neumann
algebras, integrated density of states
\vskip 0.3in
\newpage
\section{Introduction}
The goal of the Oberwolfach-Mini-Workshop ``$L^2$-spectral invariants and
the integrated density of states'' was to unify the point of views
and approaches in certain areas of geometry and mathematical
physics. The aim of our paper is to make the connection between
those fields even more explicite. Let us start with a very brief
introduction to $L^2$-spectral invariants.
\subsection{$L^2$-spectral invariants}
Let $\Gamma$ be a countable group and $L^2(\Gamma)$ be the
Hilbert-space of the formal sums $\sum_{\gamma\in\Gamma}
a_\gamma\cdot\gamma$, where $a_\gamma\in\bC$ and
$\sum_{\gamma\in\Gamma} |a_\gamma|^2<\infty$. Notice that $\Gamma$
unitarily acts on $L^2(\Gamma)$ by

\noindent
 $L_\delta(\sum_{\gamma\in\Gamma} a_\gamma\cdot\gamma)=
\sum_{\gamma\in\Gamma} a_{\gamma\delta^{-1}}\cdot\gamma\,.$ Hence
one can represent the complex group algebra $\bC\Gamma$ as bounded
operators by left convolutions. The weak closure of  $\bC\Gamma$ in
$B(L^2(\Gamma))$ is the group von Neumann algebra $\cN\Gamma$. The
group von Neumann algebra has a natural trace:
$$\tr_\Gamma(A)=\lan A(1),1\ran\,,$$
where $1\in L^2(\Gamma)$ is identified with the unit element of the
group. Let $B\in\cN\Gamma$ be a positive, self-adjoint element, then
by the spectral theorem of von Neumann
$$B=\int_{0}^{\infty}\lambda\,d E^B_\lambda\,,$$
where $E^B_\lambda=\chi_{[-\infty,\lambda]}(B)\,.$ We can associate
a {\it spectral measure} $\mu_B$ to our operator $B$ by
$$\mu_B[0,\lambda]=\sigma_B(\lambda)=\tr_\Gamma
E^B_\lambda\,.$$ Note that the jumps of $\sigma_B$ are associated to
the eigenspaces of $B$. We denote by $S^B_\lambda$ the orthogonal
projection onto the subspace of $\lambda$-eigenvectors. One can also
consider the Fuglede-Kadison determinant as
$$\det_{\Gamma}(B):=\exp(\int_0^\infty \log\lambda \,d\mu_B)\,.$$

\noindent
Now let $\Gamma$ be a finitely generated amenable group. That is
$\Gamma$ has a \Fo -exhaustion of finite
subsets $1\in F_1\subset F_2\subset\dots,
\,\cup^\infty_{n=1}F_n=\Gamma$ such that
$$\limn\frac{|KF_n|}{|F_n|}=1\,,$$
for any non-empty finite subset $K\subset \Gamma$. Then the following
approximation theorem holds:
\begin{statement} \label{s1} \cite{DLM}
For any positive, self-adjoint $B\in \bC(\Gamma)$ and $\lambda\geq
0$,
$$\tr S^B_\lambda=\limn\frac{\dim_\bC \ke (B_n-\lambda)}{|F_n|}\,,$$
where $B_n=p_{F_n} B i_{F_n}$ and $p_{F_n}:L^2(\Gamma)\to L^2(F_n)$
is the natural projection operator, $i_{F_n}:L^2(F_n)\to
L^2(\Gamma)$ is the adjoint of $p_{F_n}$, the natural imbedding
operator.
\end{statement}
As we shall see in Section \ref{unif} the fact that the
spectral distribution functions $N_{B_n}$ converge to $\sigma_B$ uniformly
immediately follows from Statement \ref{s1}.
Note that
$$N_{B_n}(\lambda)=\frac{s(\lambda,B_n)} {|F_n|}\,,$$
where $s(\lambda,B_n)$ is the number of eigenvalues of $B_n$
(counted with multiplicities) not greater than $\lambda$.

Similar approximation theorem holds for residually finite groups.
Let $\Gamma$ be a finitely generated residually finite group with
finite index normal subgroups
$$\Gamma\rhd N_1 \rhd N_2\rhd\dots, \cap^\infty_{k=1} N_k=\{1\}\,.$$
\begin{statement} \label{s2} \cite{Lueck1}
Let $B\in\bQ\Gamma$ be a positive, self-adjoint element and let
$\pi_k(B)=B_k\in\bQ(\Gamma/N_k)$ be the associated finite
dimensional linear operators, where $\pi_k:\Gamma\to \Gamma/N_k$ are
the quotient maps. Then the  distribution function $\sigma_B$ is the
uniform limit of the spectral distribution functions $N_{B_k}$.
\end{statement}
Note that the rationality assumption is crucial in the proof of
Statement \ref{s2}. For the Fuglede-Kadison determinant one has the
following result:
\begin{statement} \cite{Lueck1}
\label{s3} If $\Gamma$ is amenable or residually finite, and 

\noindent
$B\in\bQ\Gamma$ as above then $\det_\Gamma(B)\geq 1$, that is
$\int_0^\infty \log\lambda \,d\mu_B\geq 0$.
\end{statement}
\subsection{Quasi-crystal graphs}
Let $G(V,E)$ be an infinite, connected graph with bounded vertex degrees,
which is amenable in the sense that there exists a \Fo -sequence $\{F_n\}^n_{n=1}$
of finite subsets in $V$ such that
$$\limn\frac{|\partial F_n|}{|F_n|}=1\,,$$
where
$$\partial F_n:=\{p\in F_n\,\mid\,\,\mbox{ there exists $q\notin F_n$, $(p,q)\in
E$\}}\,.$$
We shall consider those graphs that exhibit a certain aperiodic order. First
we need some definitions.

\noindent
The $r$-pattern of a vertex $x\in V$ is
the graph automorphism class of the rooted ball around $x$.
That is $x,y\in V$ have the same $r$-pattern if there exists a
graph isomorphism $\phi$ between the balls $B_r(x)$ and $B_r(y)$ such
that $\phi(x)=y$. Denote by $P_r(G)$ the finite set of all possible
$r$-patterns in $G$. We say that $G$ is an {\it abstract quasicrystal graph} (AQ-graph)
if for any $\alpha\in P_r(G)$ there exists a {\it frequency} $P(\alpha)$
such that
\begin{itemize}
\item
For {\it any} \Fo -sequence $\{Q_n\}^\infty_{n=1}$:\,
$$\limn \frac{|Q^\alpha_n|}{|Q_n|}=P(\alpha)\,,$$
where $Q^\alpha_n\subseteq Q_n$ is the set of vertices in $Q_n$ having
the $r$-pattern $\alpha$. \end{itemize}
Of course, the Cayley-graphs of finitely generated groups are AQ-graphs. Also,
the dual graph of the Penrose tilings and in general, graphs defined by nice
Delone sets considered in \cite{LS1}, \cite{LS2} are AQ-graphs  without
translational symmetries.
We can construct random AQ-graphs as well. Let $G$ be the standard lattice
graph in the two-dimensional plane. Consider the Bernoulli space
$\Omega=\prod_G \{0,1\}$. For each $\omega\in \Omega$, $G_\omega$ is the
following graph. If $\omega(a,b)=0$ then we add a new edge
$((a,b),(a+1,b+1))$ to the lattice. Then for almost all $\omega\in\Omega$,
$G_\omega$ is an AQ-graph.

\noindent
The algebra that plays the role of the group algebra is the {\it algebra
of pattern invariant matrices} $P_G$. A pattern-invariant matrix is a function
$A:V\times V\to\bC$ satisfying the following properties:
\begin{itemize}
\item
$A(x,y)=0$ if $d_G(x,y)>r_{A}$, where $d_G$ is the shortest path
distance on the vertices. The value $r_{A}$, the propagation of $A$,
depends only on $A$.
\item
$A(x,y)=A(\phi(x),\phi(y))$ if $\phi$ is a graph isomorphism (there
can be more than one) between the $r_{\cA}$-ball around $x$ and the
$r_{A}$-ball around $\phi(x)$, mapping  $x$ to $\phi(x)$.
\end{itemize}
Of course, \begin{itemize}
\item $AB(x,y)=\sum_{z\in V}A(x,z)B(z,y)$  \item  $ (A+B)(x,y)=A(x,y)+B(x,y)$ \item
$(\lambda A)(x,y)=\lambda(A(x,y))$ \item  $A^*(x,y)=\overline{A(y,x)}$
\end{itemize} defines a
  $\star$-algebra
structure on $P_G$.
Note that $P_G$ is naturally represented on $L^2(V)$ by
$$A(f)(x)=\sum_{y\in V} A(x,y) f(y)\,,$$ for $f\in L^2(V)$.
These sort of operators were considered e.g. in \cite{Elek},
\cite{LS1},\cite{LS2}.

\noindent In the case of tilings one can imbed $P_G$ into a von
Neumann algebra defined by the associated Delone system. In
\cite{Elek}, we imbedded $P_G$ into a continuous ring defined by
the special construction of certain self-similar graphs. In this
paper we define a von Neumann algebra $N_G$ using the natural trace
on $P_G$ and consider the imbedding  $\psi_G:P_G\to N_G$ (see
Section \ref{pattern}). Now we can state the main results of this
paper.
\begin{theorem} \label{t1}
Let $G$ be an AQ-graph, $B\in P_G$ a positive, self-adjoint operator
(as an operator on $L^2(V)$. Suppose that for any $x,y\in V$,\,$
B(x,y)\in \bQ$. Then for any \Fo -sequence $\{Q_n\}^\infty_{n=1}$,
the spectral distribution functions $N_{B_n}$ uniformly converge to
the integrated density of states (IDS) function $\sigma_B$, where
$B_n=p_{Q_n} B i_{Q_n}$ and $\sigma_B$ is the von Neumann spectral
function of the positive, self-adjoint element $\psi_G(B)$.
\end{theorem}

 Note that in  \cite{LS1},
\cite{LS2} similar uniform convergence results was proved for a
class of AQ-graphs without the rationality assumption. In
\cite{Elek}, we also considered certain AQ-graphs for which uniform
convergence follows without the rationality assumption.
\begin{theorem} \label{t2}
If $G$ and $B$ are as above then
$$\int^\infty_0 \log \lambda \,d\mu_B \geq 0\,,$$
where $\mu_B$ is the spectral measure associated to $\psi_G(B)$.
\end{theorem}
As a consequence of this theorem we have the following corollary:
\begin{corollary}
\label{cor} If $G$ and $B$ are as above and $\lambda$ is a point of
discontinuity of the IDS function $\sigma_B$, then $\lambda$ is an
algebraic number.
\end{corollary}

\section{Uniform spectral convergence}\label{unif}
In this section we prove some technicalities on the convergence of spectral
distribution functions.
\begin{proposition}
\label{p10}
Let $\{f_n\}^\infty_{n=1},f$ be monotone increasing right-continuous functions
on the interval $[0,K]$ such that
$f_n(K)=1$ for any $n\geq 1$ and $f(K)=1$ as well. Suppose
that $\limn f_n(\lambda)=f(\lambda)$ if $\lambda$ is a point of continuity of
$f$.
Also suppose that if $\lambda$ is a point of discontinuity of $f$, then
$$\limn J_n(\lambda)=J(\lambda)\,,$$
where $$J_n(\lambda)=f_n(\lambda)-f^-_n(\lambda)\quad
J(\lambda)=f(\lambda)-f^-(\lambda)$$
are the jumps of our functions.
Then $\{f_n\}^\infty_{n=1}$ uniformly converge to $f$.
\end{proposition}
\proof
First we need the following lemma.
\begin{lemma}
\label{l11}
$$\limn f^-_n(y)=f^-(y)$$ for any $y\in [0,K]$.
\end{lemma}
\proof
First suppose that
\begin{equation}
\label{e11}
f^-(y)-\liminf_{n\to\infty} f^-_n(y)=\delta>0\,.
\end{equation}
Let $z<y$ be a continuity point of
 $f$ such that $f^-(y)-f(z)<\frac{\delta}{2}\,.$
Then for large enough $n$, $|f_n(z)-f(z)|<\frac{\delta}{2}.$ Since
$f_n(z)\leq f_n^-(y)$ (\ref{e11}) leads to a contradiction.
Now suppose that
\begin{equation}
\label{e11b}
\limsup_{n\to\infty} f^-_n(y)-f^-(y) =\delta>0\,.
\end{equation}
Let $y<z$ a continuity point of $f$ such that
$f(z)-f(y)<\frac{\delta}{10}\,.$ Again, for large enough $n$ ;
$$|f_n(z)-f(z)|< \frac{\delta}{10} \,\,\,\mbox{and}\,\,\,
|(f_n(y)-f_n^-(y))- (f(y)-f^-(y))|<\frac{\delta}{10}\,.$$
Since $f_n(z)\geq f_n(y)$, (\ref{e11b}) leads to a contradiction. \qed

\noindent
Now we turn to the proof of our proposition. If the uniform convergence
does not exist, then there are points $\{x_{n_k}\}^\infty_{k=1}, y\in [0,K]$
such that $x_{n_k}\to y$ and
\begin{equation}
\label{e12}
|f_{n_k}(x_{n_k})-f(x_{n_k})|\geq \e >0\,
\end{equation}
 Let $z<y$ be a continuity point of
 $f$ such that $f^-(y)-f(z)<\frac{\e}{10}\,.$
If $k$ is large enough, then
$$|f_{n_k}(z)-f(z)|<\frac{\e}{10}
\,\,\,\mbox{and}\,\,\,
|f^-_{n_k}(y)-f^-(y)|<\frac{\e}{10}\,.$$
Hence for large enough $k$:
$$|f_{n_k}(x_{n_k})-f(x_{n_k})|< \e\,,$$
if $z\leq x_{n_k} \leq y$.
Now let $y<z$ be a continuity point of
 $f$ such that $f(z)-f(y)< \frac{\e}{10}\,.$
Again if $k$ is large, then
$$|f_{n_k}(z)-f(z)|<\frac{\e}{10}\,\,,
|f^-_{n_k}(y)-f^-(y)|<\frac{\e}{10}\,\,\mbox{and}\,\,
|J_{n_k}(y)-J_{n_k}(y)|< \frac{\e}{10}\,.$$
Hence for large $k$, $|f_{n_k}(x_{n_k})-f(x_{n_k})|<\e$ if $y\leq
x_{n_k}\leq z$. Hence (\ref{e12}) leads to a contradiction. \qed

\noindent
Suppose that $\mu_n,\mu$ are probability measures on the interval $[0,K]$,
and $\{\mu_n\}^\infty_{n=1}$ weakly converge to $\mu$. Then if
$F_n(\lambda)=\mu_n([0,\lambda]$, $F(\lambda)=\mu([0,\lambda]$, then
$\limn F_n(\lambda)=F(\lambda)$ if $\lambda$ is a continuity point of $F$
that is when $\mu(\{\lambda\})=0\,.$ Suppose that
\begin{itemize}
\item $\mu_n$ weakly converge to $\mu$ and
\item if $\mu(\{\lambda\})>0$, then $\mu_n(\{\lambda\})\to \mu (\{\lambda\})$.
\end{itemize}
Then by our Proposition, $\{F_n\}^\infty_{n=1}$ converge to $F$ uniformly.

\noindent
\underline{\bf Remark:} \, Consider the situation in Statement \ref{s1}.
If $B\in\bC\Gamma$ is a positive, self-adjoint element, then by \cite{DLM}
$$\tr_\Gamma S^B_\lambda=\limn \frac{\dim_{\bC} \ke(B-\lambda)}{|F_n|}\,,$$
where $S^B_\lambda$ is the orthogonal projection onto the $\lambda$-eigenspace
of $B$. Hence by Proposition \ref{p10}, $N_{B_n}$ uniformly converge to the
spectral measure $\sigma_B$ (the pointwise convergence was considered
in \cite{Schick}). 
\section{The pattern-frequency algebra construction} \label{pattern}
Let $G$ be an AQ-graph and $P_G$ be the $\star$-algebra of pattern-invariant
matrices. Then for $A\in P_G$ we define
\begin{equation} \label{e15}
\tr_G(A):=\lim\frac{1}{|Q_n|}\sum_{x\in Q_n} A(x,x)\,,
\end{equation}
where $\{Q_n\}\i$ is a \Fo -sequence in $G$.
\begin{proposition}
\label{p15}
\begin{description}
\item[(a)] The limit in (\ref{e15}) exists.
\item[(b)] $\tr_G(A)$ does not depend on the choice of the \Fo -sequence.
\item[(c)] $\tr_G$ is a trace, that is a linear functional satisfying
 $\tr_G(AB)=\tr_G(BA)$.
\item[(d)] $\tr_G$ is faithful, that is $\tr_G(A^*A)>0$ if $A\neq 0$.
\end{description}
\end{proposition}
\proof
Note that
$$\frac{1}{|Q_n|}\sum_{x\in Q_n} A(x,x)= \frac{1}{|Q_n|}
\sum_{\alpha\in P_{r_A}(G)}
|Q_n^\alpha| A_\alpha\,,$$
where $Q_n^\alpha$ is the number of vertices in $Q_n$ with $r_A$-pattern
$\alpha$ and $A_\alpha=A(x,x)$ if $x$ has $r_A$-pattern $\alpha$. Since
$\limn \frac{|Q_n^\alpha|}{|Q_n|}=P(\alpha)$, (a) and (b) immediately follows.

\noindent
Now let $A,B\in P_G$, then
$$\tr_G(AB)=\limn (\frac{1}{|Q_n|}\sum_{x\in Q_n}\sum_{y\in V} A(x,y)
B(y,x))$$
$$\tr_G(BA)=\limn (\frac{1}{|Q_n|}\sum_{x\in Q_n}\sum_{y\in V} B(x,y)
A(y,x))$$
That is \begin{equation}\label{e16}
|\tr_G(AB)-\tr_G(BA)|\leq \limn \frac{2}{|Q_n|}  \sum_{x\in\partial_D Q_n}
m_A m_B\,,
\end{equation}
where $m_A=\sup_{x,y} |A(x,y)|, m_B=\sup_{x,y} |B(x,y)|\,,$ and
$$\partial_D Q_n:=\{x\in Q_n\,\mid \mbox{there exists}\, y\notin Q_n,
d_G(x,y)\leq \max(r_A,r_B)\}\,.$$
Note that $m_A,m_B<\infty$ by the pattern-invariance.
Also, $\limn \frac {|\partial_D Q_n|}{|Q_n|}=0\,$ by the \Fo -property and the
bounded
vertex degree condition.
Hence by (\ref{e16}) (c) follows. Finally, let $0\neq A\in P_G$.
Then
$$\tr_G(A^*A)=\frac{1}{|Q_n|}\sum_{x\in Q_n}(\sum_{y\in V} A(x,y)
\overline{A(y,x)})= \frac{1}{|Q_n|}\sum_{x\in Q_n}(\sum_{y\in V}
|A(x,y)|^2)\,.$$
That is $\tr_G(A^*A)=\sum_{\alpha\in P_{r_A}(G)} P(\alpha) L_\alpha\,,$
where $L_\alpha=\sum_{y\in V}
|A(x,y)|^2$ if $x$ has the $r_A$-pattern $\alpha$. Therefore (d) follows. \qed

\noindent
Now the pattern-frequency algebra is constructed by the  GNS-construction the
usual way. The algebra $P_G$ is a pre-Hilbert space with respect to
the inner product $<A,B>=\tr_G(B^*A)\,.$
Then $L_A B=AB$ defines a representation of $P_G$ on this pre-Hilbert space.
\begin{lemma}
\label{l18}
$L_A$ is a bounded operator if $A\in P_G$. \end{lemma}
\proof
Let $A,B\in P_G$. Denote by $\|\,\|$ the pre-Hilbert space norm. Then
$$\|B\|^2=\limn (\frac{1}{|Q_n|}\sum_{x\in Q_n}(\sum_{y\in V}|B(x,y)|^2)\,.$$
$$\|L_A B\|^2=\limn (\frac{1}{|Q_n|}\sum_{x\in Q_n}(\sum_{y\in V}
BB^*(x,y)A^*A(y,x))\leq \limn (\frac{K}{|Q_n|}\sum_{x\in Q_n}(\sum_{y\in V}
|BB^*(x,y)|)\,,$$
where $K=\sup_{x,y} |A^*A(x,y)|$.
Note that
$$|BB^*(x,y)|=|\sum_{z\in V} B(x,z) B^*(z,y)|\leq
\sum_{z\in V}|B(x,z)| |B(y,z)|\leq $$ $$ \leq\frac{1}{2}
\sum_{z\in V}(|B(x,z)|^2+|B(y,z)|^2)\,.$$
Let $t_x:=\sum_{y\in V} |B(x,y)|^2\,.$ Then
$$\|B\|^2=\limn \frac{1}{|Q_n|}\sum_{x\in Q_n}t_x\,.$$
$$\|L_A B\|^2=\limn (\frac{1}{|Q_n|}\sum_{x\in Q_n}(\sum_{y\in V}\frac{1}{2}
(t_x+t_y)\leq \limn \frac{KM} {|Q_n|}\sum_{x\in Q_n}t_x\,,$$
where $M$ is the maximal number of vertices in an $r_{A^*A}$-ball of
$G$. Hence $\|L_A B\|^2\leq KM \|B\|^2\,.$\qed

\noindent
Thus $P_G$ is represented by bounded operators on the Hilbert-space
closure. Now the pattern-frequency algebra $N_G$ is defined as the
weak-closure of $P_G$. then $\tr_G(A)=\langle L_A 1, 1 \rangle$ extends to
$N_G$ as an ultraweakly continuous, faithful trace. 
The finite von Neumann algebra $N_G$ shall play the
role of the group von Neumann algebra in our paper.

\section{Spectral approximation}
\begin{proposition}
\label{p21} If $B$ is a positive, self-adjoint operator in $P_G$,
then $\psi_G(B)$ is positive self-adjoint as an element of $N_G$.
\end{proposition} \proof The self-adjointness is obvious, since the
representation $\psi_G:P_G\to N_G$ preserves the $*$-structure:
$$\langle L_B X,Y\rangle=\tr_G(Y^* B X)\,$$
$$\langle X, L_B Y \rangle=\tr_G((BY)^*X)=\tr_G(Y^* B X)\,.$$
The operator $B$ is positive in $P_G$ if and only if
$\langle B(f),f\rangle_{L^2(V)}\geq 0$ for any $f\in L^2(V)$. That is
\begin{equation}
\label{e21}
\sum_{y,z\in V} B(y,z) f(z) \overline{f(y)}\geq 0\,.
\end{equation}
On the other hand, $\psi_G(B)$ is positive if
$\tr_G(U^*BU)\geq 0$ for any $U\in P_G$.
$$\tr_G(U^*BU)=\limn (\frac{1}{|Q_n|}\sum_{x\in Q_n}\sum_{y,z\in
  V}U^*(x,y)B(y,z)U(z,x)= $$ $$=
\limn (\frac{1}{|Q_n|}\sum_{x\in Q_n}\sum_{y,z\in V} B(y,z) U(z,y)
\overline{U(y,x)})\,.$$
By (\ref{e21}), for any $x\in Q_n$:
$$\sum_{y,z\in V} B(y,z) U(z,y)
\overline{U(y,x)}\geq 0\,.$$
Hence the positivity follows. \qed

\noindent
Now we need some estimates for the norms of the finite dimensional
linear transforms $B_n$.
\begin{proposition}
\label{p22}
If $B\in P_G$ is positive, self-adjoint operator then
$B_n$ is positive, self-adjoint as well. Furthermore, there exists $K>0$ such
that $\|B_n\|\leq K$, for any $n\geq 1$, where the norm of $B_n$ is considered
as a linear transformation in $\mat_{Q_n\times Q_n}(\bC)\,.$ \end{proposition}
\proof
Recall that $B_n=p_{Q_n}B p^*_{Q_n}$, hence the
positivity and self-adjointness are obvious. The norm-estimate immediately
follows from the next lemma.
\begin{lemma}
Let $H(V,E)$ be a finite graph, where the vertex degrees are not greater than
$d$. Suppose that $A:V\times V\to \bC$ is a matrix such that
\begin{itemize}
\item $|A(x,y)|\leq m$
\item $A(x,y)=0$ if $d_H(x,y)>r$. \end{itemize}
 Then $\|A\|\leq C_{d,m,r}$, where the
constant $C_{d,m,r}>0$ depends only on the parameters, not on the graph $H$
itself.
\end{lemma}
\proof
For $f,g\in L^2(V)$ unit vectors,
$$|\langle A(f),g\rangle|=|\sum_{x,y\in V} A(x,y) f(x) \overline{g(y)}|\leq
m |\sum_{x,y\in V\,,d_H(x,y)\leq r}  f(x)\overline{g(y)}|\leq $$ $$\leq
m \sum_{x,y\in V\,,d_H(x,y)\leq r}  ( |f(x)|^2+|g(y)|^2)\,.$$
The number of occurences of $|g(y)|^2$ on the right hand side is not greater
than the maximal possible number $T(r,d)$ of vertices in a ball of radius $r$
in a graph of vertex degrees bounded by $d$. This quantity depends only on $r$
and $d$.
Thus
$$\|A\|=\sup_{f,g\in L^2(V), \|f\|=1, \|g\|=1} |\langle A(f),g\rangle|\leq m
(1+ T(r,d))\,. \quad\qed.$$

\noindent
Now we prove our key trace approximation formula.
\begin{proposition}
\label{p25}
Let $B\in P_G$, then for any $k\geq 1$,
$$\tr_G(B^k)=\limn \frac{1}{|Q_n|} \tr(B_n^k)\,.$$ \end{proposition}
\proof
$$\tr_G(B^k)= \limn \frac{1}{|Q_n|}\sum_{x\in Q_n} B^k(x,x)= $$ $$=
\limn \frac{1}{|Q_n|}\sum_{x\in Q_n} \sum_{y_1,y_2,\dots,y_{k-1}\in V}
B(x,y_1) B(y_1,y_2)\dots B(y_{k-1},x)\,.$$
On the other hand,
$$\frac{1}{|Q_n|} \tr(B_n^k)= \frac{1}{|Q_n|}\sum_{x\in Q_n}\sum_{y_1,y_2,\dots,y_{k-1}\in Q_n}
B(x,y_1) B(y_1,y_2)\dots B(y_{k-1},x)\,.$$
Note that if $x,y\in Q_n$ then $B(x,y)=B_n(x,y)$.
Clearly,
$$|\frac{1}{|Q_n|}\sum_{x\in Q_n} \sum_{y_1,y_2,\dots,y_{k-1}\in V}
B(x,y_1) B(y_1,y_2)\dots B(y_{k-1},x)-
$$  $$-\frac{1}{|Q_n|}\sum_{x\in Q_n} \sum_{y_1,y_2,\dots,y_{k-1}\in Q_n}
B(x,y_1) B(y_1,y_2)\dots B(y_{k-1},x)|\leq \frac{1}{|Q_n|} m_B^k s_n\,,$$
where $m_B=\sup_{x,y\in V} B(x,y)$ and $S_n$ is the set of strings
$\{x,y_1,y_2,\dots,y_{k-1}\}$ such that
\begin{itemize}
\item $d_G(x,y_1)\leq r$\,.
\item $d_G(y_i,y_{i+1})\leq r$ for any $1\leq i \leq k-2$.
\item For some $1\leq i \leq k-1$: $y_i\notin Q_n$. \end{itemize}
Therefore our proposition follows from the lemma below.
\begin{lemma}
\label{l26}
$\limn \frac{|S_n|}{|Q_n|}=0\,.$
\end{lemma}
\proof
For $a\in \bN$, let again
$$\partial_a Q_n:=\{x\in Q_n\,\mid \mbox{there exists}\, y\notin Q_n,
d_G(x,y)\leq a \}\,.$$
By the \Fo -property and the vertex degree condition $\limn \frac{|\partial_a
  Q_n|}{|Q_n|}=0\,.$
If $x$ is a starting vertex of a string in $S_n$, then $x\in\partial_{kr}
Q_n$. The number of choices for $y_1$ is not greater than $q_r$, where $q_r$
is the maximal number of vertices in an $r$-ball of $G$. Hence
$$|S_n|\leq |\partial_{kr} Q_n| |q_r|^{k-1}\,.$$
Thus our lemma follows. \qed
\section{The integrated density of states}
Let $G$ be an AQ-graph, $B\in P_G$ be positive, self-adjoint and
$\{B_n\}\i$ be the associated, finite dimensional linear transformations.
Since $\psi_G(B)\in N_G$ is positive, self-adjoint element of
a von Neumann-algebra we have the spectral theorem:
$$\psi_G(B):=\int_0^\infty\lambda \,dE^B_\lambda\,,$$
where $E^B_\lambda$ are the associated spectral projections
$\chi_{[0,\lambda]}(\psi_G(B))\,.$ The spectral measure is
defined by
$$\sigma_B(\lambda)=\mu_B[0,\lambda]:=\tr_G(E^B_\lambda)\,.$$
For the finite dimensional operators $B_n$
$$N_{B_n}(\lambda)=\frac{1}{|Q_n|} \tr(E^{B_n}_\lambda)\,$$
are the usual spectral distributions.
Let $K>0$ be a real number such that $\|B\|,\|B_n\|<K$.
Then if $P\in \bR[x]$ is a polynomial,
$$\int^K_0 P(\lambda) d\mu_{B_n}(\lambda)=\frac{1} {|Q_n|} \tr (P(B_n))\,.$$
$$\int^K_0 P(\lambda) d\mu_{B}(\lambda)=\tr_G ( P(B_n))\,.$$
Hence by Proposition \ref{p25}
$$\limn \int^K_0 P(\lambda) d\mu_{B_n}(\lambda)=
\int^K_0 P(\lambda) d\mu_{B}(\lambda)\,.$$
That is the probability measures $\mu_{B_n}$ weakly converge to
$\mu_B$. In other words, $\sigma_B$ is the integrated density of states.
This immediately implies the associated Shubin-formula:
$$\limn N_{B_n}(\lambda)=\tr_G E^B_\lambda\,,$$
for any $\lambda\geq 0$ continuity point of $\sigma_B$.
\section{Approximation of the ground state density}
The goal of this section is to prove the following proposition.
\begin{proposition}
\label{p33} If $B$ is a positive, self-adjoint element of $P_G$
 and $B(x,y)$ is a rational number for any
$x,y\in V$, then
$$\limn \frac{\dim_{\bC} \ke B_n}{|Q_n|}=\tr_G E^B_0\,.$$
\end{proposition}
\proof The proof is very similar to Theorem 6.9 \cite{Schick}
nevertheless there are
some details we would like to make precise for the reader who is
not familiar with original paper of L\"uck \cite{Lueck1}.
Note that $B$ has only finitely many different values, therefore we can
suppose (after multiplying by an appropriate integer) that $B$ has only
integer values.

First we need some notations.
For a monotone function $f$,
$$f^+(\lambda)=\inf_{\e\to 0} f(\lambda+\e)\,$$
$$\osb(\lambda):=\limsup_{n\to\infty} N_{B_n}(\lambda)$$
$$\usb(\lambda):=\liminf_{n\to\infty} N_{B_n}(\lambda)\,.$$
Let us recall Lemma 6.5 \cite{Schick} in the special case we need it.
\begin{lemma}
\label{l30}
For any polynomial $P_k\in \bR[x]$ which has the property
$\chi_{[0,\lambda]}(x)\leq P_k(x)\leq \frac{1}{k} \chi_{[0,\lambda]}(x)+
\chi_{[0,\lambda+\frac{1}{k}]}(x)\,$ for any $0\leq x \leq K$,
we have the following inequalities:
$$\sigma_B(\lambda)\leq \tr_G (P_k(\psi_G(B))\leq
\sigma_B(\lambda+\frac{1}{k})+\frac{1}{k}\,.$$
$$N_{B_n}(\lambda)\leq \frac{1}{|Q_n|} \tr(P_k(B_n))\leq
N_{B_n}(\lambda+\frac{1}{k})+\frac{1}{k}\,.$$ \end{lemma}
We also need the appropriate version of Proposition 6.7 \cite{Schick}.
\begin{proposition}
\label{p31}
$$\osb(\lambda)\leq \sigma_B(\lambda)=\usb^+(\lambda)=\osb^+(\lambda)\,.$$
\end{proposition}
\proof
First choose the polynomial $P_k$ as in the previous lemma. Then
$$N_{B_n}(\lambda)\leq \frac{1}{|Q_n|} \tr(P_k(B_n))\leq
N_{B_n}(\lambda+\frac{1}{k})+\frac{1}{k}\,.$$
That is as $n\to\infty$:
$$\osb(\lambda)\leq \tr_G P_k(\psi_G(B))\leq \usb (\lambda+\frac{1}{k})+
\frac{1}{k}\,.$$
Now we let $k\to\infty$. Then:
$$\lim_{k\to\infty} \tr_G P_k(\psi_G(B))=
\lim_{k\to\infty} \langle P_k(\psi_G(B)) (1),1\rangle=
\langle \chi_{[0,\lambda]} (\psi_G(B)) (1),1\rangle=\tr_G E^B_\lambda\,,$$
since $\{P_k(\psi_G(B))\}^\infty_{k=1}$ strongly converge to
$\chi_{[0,\lambda]}(\psi_G(B)\,.$
Therefore,
$$\osb(\lambda)\leq \sigma_B(\lambda)\leq \usb^+(\lambda)\,.$$
Clearly,
$$\sigma_B(\lambda)\leq \usb(\lambda+\e)\leq \osb(\lambda+\e)\leq
\sigma_B(\lambda+\e)\,.$$
That is
$$\sigma_B(\lambda)=\usb^+(\lambda)=\osb^+(\lambda)\,.$$

Now we need a lemma on Riemann-Stieltjes integral. The proof
can be found in \cite{Lueck2} in  Lemma 4.1.
\begin{lemma} \label{stieltjes}
Let $f$ be a continuously differentiable function on the positive reals
and $\mu$ be a probability measure on the $[0,K]$ interval. Suppose
that $F$ is the distribution function of $\mu$, that is $\mu[0,\lambda]=
F(\lambda)$. Then for any $0<\e\leq K$:
$$\int_\e^K f(\lambda) d\mu= -\int^K_\e f'(\lambda)F(\lambda) d\lambda
+ f(K)F(K)-f(\e)F(\e)\,.$$
\end{lemma}
Let $|det|_1(B_n)$ is defined as the product of non-negative
eigenvalues of $B_n$, this is the lowest nonzero coefficient in the
characteristic polynomial of $B_n$ hence a non-zero integer.
Therefore,
\begin{equation}
\label{logdet}
\log \det_{Q_n} (B_n):=\frac{1}{|Q_n|}\log (|det|_1(B_n)\geq 0.
\end{equation}
By Lemma \ref{stieltjes},
$$\log \det_{Q_n} (B_n)=\log K- \log \e \cdot N_{B_n}(\e)-
\int^K_\e N_{B_n}(\lambda) \frac{1}{\lambda} d\lambda\,,$$
where $\e>0$ is a real number below the smallest eigenvalue of $B_n$.

\noindent
Since $\int^K_\e \frac{N_{B_n}(0)}{\lambda} d\lambda=(\log K -\log \e)
N_{B_n}(0)\,$:
\begin{equation}
\label{e33a}
\log \det_{Q_n} (B_n)=\log K(1-N_{B_n}(0))-\int^K_{0}
\frac{N_{B_n}(\lambda)- N_{B_n}(0)}{\lambda} d\lambda.
\end{equation}
Hence by (\ref{logdet})
\begin{equation}
\label{e33b}
\int^K_{0}
\frac{N_{B_n}(\lambda)- N_{B_n}(0)}{\lambda} d\lambda\leq \log K\,.
\end{equation}
Now we estimate $\int^K_{0}
\frac{\sigma_{B}(\lambda)- \sigma_{B}(0)}{\lambda} d\lambda$.
$$
\int^K_{\e}
\frac{\sigma_{B}(\lambda)- \sigma_{B}(0)}{\lambda} d\lambda\leq
\int^K_{\e}\frac{\usb^+(\lambda)- \sigma_{B}(0)}{\lambda} d\lambda\leq
\int^K_{\e}\frac{\usb(\lambda)- \sigma_{B}(0)}{\lambda} d\lambda\,.$$
Furthermore,
$$
\int^K_{\e}\frac{\usb(\lambda)- \sigma_{B}(0)}{\lambda} d\lambda\,\leq
\int^K_{\e}\frac{\usb(\lambda)- \osb(0)}{\lambda} d\lambda=
\int^K_{\e}\frac{\liminf_{n\to \infty} N_{B_n}(\lambda)-
N_{B_n}(0)}{\lambda} d\lambda\,\leq$$
$$\leq \liminf_{n\to \infty} \int^K_\e \frac{ N_{B_n}(\lambda)-
N_{B_n}(0)} {\lambda} d\lambda\leq \log K\,.$$
Now we take the limit $\e\to 0$:
$$
\int^K_{0}
\frac{\sigma_{B}(\lambda)- \sigma_{B}(0)}{\lambda} d\lambda\leq$$
\begin{equation} \label{e34} \leq
\int^K_{0}
\frac{\usb(\lambda)- \osb(0)}{\lambda} d\lambda\leq
\sup_{\e\to 0} \liminf _{n\to \infty} \frac {N_{B_n}(\lambda)-
N_{B_n}(0)}{\lambda} d\lambda\leq \log K\,.
\end{equation}
Therefore $
\lim_{\lambda\to 0} \usb(\lambda)=\osb(0)\,.$
That is by Proposition \ref{p31}:
$\osb(0)=\sigma_B(0)$, hence \begin{equation} \label{this}
\limsup_{n\to \infty} N_{B_n}(0)
=\sigma_B(0)\,.\end{equation}
Since one can consider any subsequence of $B_n$, (\ref{this}) implies
that
$$\limn N_{B_n}(0)=\limn\frac{\dim_{\bC} \ke B_n}{|Q_n|}=\sigma_B(0)
=\tr_G E^B_0\,.\quad\qed$$
\section{The proof of Theorem \ref{t2}}
We would like to show that
$\int^\infty_0 \log \lambda d\mu_B\geq 0\,.$
Note that
$$\int^\infty_0 \log \lambda d\mu_B=
\lim_{\e\to 0} \int^K_{\e} \log \lambda d\mu_B=
\lim_{\e\to 0} (\log K (1-\sigma_B(\e))-\int^K_{\e}
\frac{\sigma_B(\lambda)-\sigma_B(0)}{\lambda} d\lambda\,.$$
Hence
$$\int^\infty_0 \log \lambda d\mu_B=\log K (1-\sigma_B(0))-
\int^K_{0}
\frac{\sigma_B(\lambda)-\sigma_B(0)}{\lambda} d\lambda\,.$$
By (\ref{e34})
$$\sup_{\e>0}\liminf_{n\to\infty} \int^K_{\e}
\frac {N_{B_n}(\lambda)-
N_{B_n}(0)}{\lambda} d\lambda \leq
\liminf_{n\to\infty}\sup_{\e>0}\int^K_{\e}
\frac {N_{B_n}(\lambda)-
N_{B_n}(0)}{\lambda} d\lambda=$$
$$=\liminf_{n\to\infty} \int^K_{0} \frac {N_{B_n}(\lambda)-
N_{B_n}(0)}{\lambda} d\lambda\leq \log K (1- N_{B_n}(0))\,.$$
Thus,
$$\log \det (B)= \log K (1-\sigma_B(0))-\int^K_0
\frac{\sigma_B(\lambda)-\sigma_B(0)}{\lambda} d\lambda\geq $$ $$
\log k (1 - \limn N_{B_n}(0))-\liminf_{n\to\infty}
\int^K_{0} \frac {N_{B_n}(\lambda)-
N_{B_n}(0)}{\lambda} d\lambda=$$
$$=\limsup_{n\to \infty} (\log (K) (1 - \limn N_{B_n}(0))-
\int^K_{0} \frac {N_{B_n}(\lambda)-
N_{B_n}(0)}{\lambda} d\lambda)= $$ $$=\limsup_{n\to \infty}
\log \det_{Q_n} (B_n)\geq 0\,. \quad \qed$$
\section{The proof of Theorem \ref{t1}}
By Proposition \ref{p10}, it is enough to prove that
$$\limn \frac{\dim_{\bC}\ke (B_n-\lambda)}{|Q_n|}=\tr_G S^B_\lambda\,,$$
where $S^B_\lambda$ is the orthogonal projection onto the $\lambda$-eigenspace
of $\psi_G(B)$. Observe that $S^B_\lambda= S_0^{(B-\lambda)^2}$. Therefore
by Proposition \ref{p33}, Theorem \ref{t1} follows from the lemma below.
\begin{lemma}
$$\limn \frac{\dim_{\bC} \ke (B_n-\lambda)-\dim_{\bC} \ke
(p_{Q_n}(B_\lambda)^2i_{Q_n})} {|Q_n|}=0\,.$$
\end{lemma}
\proof
Let $Q_n'=Q_n\backslash \partial_{2r_{B_n}} Q_n\,.$
Observe that if $\supp f\subseteq Q_n',$ then
$$(p_{Q_n}(B_\lambda)^2i_{Q_n}) f= (p_{Q_n}(B_\lambda)i_{Q_n})^2\,.$$
Hence
$$\ke (B_n-\lambda)\cap L^2(Q_n')=
\ke (p_{Q_n}(B_\lambda)^2i_{Q_n})\cap  L^2(Q_n')\,.$$ Since $\limn
\frac{|Q_n'|}{|Q_n|}=1\,$ the lemma follows. \qed

\end{document}